\documentclass{article}
\usepackage{amsmath,amssymb,latexsym}

\newtheorem{thm}{Theorem}[section]
\newtheorem{lem}{Lemma}[section]

\newcommand{\mc}{\mathcal}
\newcommand{\mb}{\mathbb}

\begin{document}

\title{On the spectrum of a differential operator on a Hilbert-P\'olya space}

\author{Xian-Jin Li\\
Department of Mathematics\\
Brigham Young University\\
Provo, Utah 84602 USA\\
xianjin@byu.edu}
\date{ }
\maketitle

\begin{abstract} In this paper we study the spectrum of a fundamental differential 
operator on a Hilbert-P\'olya space.  A number is an eigenvalue of this differential 
operator if and only if it is a nontrivial zero of the Riemann zeta function.  
An explicit formula is given for the eigenfunction associated with each nontrivial 
zero of the zeta function.  Every eigenfunction is characterized via the Poisson 
summation formula by a sequence of mysterious functions
whose explicit formulas are given.
    \end{abstract}

 \noindent{\bf Subject Class:} Primary 30A99


  \section{Introduction}

     Let $S(\mb R_+^\times)$ be the set of all infinitely
differentiable functions $f$ on $\mb R_+^\times$
with $\sup_{x\in(0,\infty)}|(\log x)^a f^{(b)}(x)|<\infty$ for $a,b=0,1,2,\cdots$.
The fundamental differential operator on $S(\mb R_+^\times)$ is
 \[Df(x)=-xf^\prime(x).\]
 We denote
  $Zf(x)=\sum_{n=1}^\infty f(nx)$,
   \[\mc H_-=\{f\mid\mb R_+^\times\to\mb C
 \text{ such that }  x^\alpha f(x)\in S(\mb R_+^\times)
  \text{ for every real number $\alpha$}\},\]
  and $\mc H_\cap$ is the set of all even functions
 $f\in S(\mb R)$ such that $f(0)=0$ and $\frak Ff(0)=0$
 where the Fourier transform of $f$ is
  \[\frak F f(y)=\int_{-\infty}^\infty f(x)e^{2\pi ixy}dx.\]
 The the quotient space
  \[\mc H=\mc H_-/Z\mc H_\cap\] is called a Hilbert-P\'olya space.
  Let $D_-$ be the operator induced by $D$ on $\mc H$.

   If $\rho$ is a nontrivial zero of $\zeta(s)$, the author 
   proved in \cite{li3} that
 $$F_\rho(x)=\int_1^\infty Z\eta(tx)t^{\rho-1}dt$$
  is an eigenfunction of $D_-$
 in $\mc H$ associated with the eigenvalue
 $\rho$, where $\eta(x)=8\pi x^2\left(\pi x^2-\frac 32\right)e^{-\pi x^2}$.
 That is,
    \begin{equation}\label{eq1.1}
    DF_\rho(x)=\rho F_\rho(x)+Z\eta(x),\quad
  \eta\in\mc H_\cap.
  \end{equation}
   It is also proved in \cite{glwx} that a complex number $
   \rho$ is an eigenvalue of $D_-$ on $\mathcal H$ if and only if $\rho$ is
  a nontrivial zero of the zeta function.

  In his research on global trace formula \cite{li2}, the author noticed a mysterious
  function associated with each function in $\mc H$ that
  prevents the Poisson summation formula to hold for that function.  In \cite{li3}
  he was able to identify this mysterious function and obtained the following result.

  \begin{thm}\label{thm1.1} (\cite[Theorem 1.3]{li3})  Let
  $\delta_{l,F}(x)=\int_0^\infty Z^{-1}F(xt)\, \frac{\sin(2l+1)\pi t}{\pi t}dt$. Then
 \[\lim_{l\to\infty}\left(\frac 1x\, \sum_{m=1}^l
 \frak FZ^{-1}F\left(\frac mx\right) -\delta_{l,F}(x)\right)=F(x)\]
for every $F\in\mc H_-$.
 \end{thm}

    Because of the appearance of the unbounded operator $Z^{-1}$
    inside $\delta_{l,F}(x)$, it is difficult to know what the function
    $\delta_{l,F}(x)$ is exactly.  In this paper, author was able to get an
    exact simple formula for $\delta_{l,F_\rho}(x)$ for every
    nontrivial zero $\rho$ of the Riemann zeta function.
     More precisely, the author obtained the following formula.

   \begin{thm} \label{thm1.2} For every nontrivial zero $\rho$ of the
   Riemann zeta function, we have
       \begin{align*}\delta_{l,F_\rho}(x)
       &={\widehat\eta(1-\rho)(l+{1\over 2})^\rho
   \over \rho x^\rho}+4\pi ({l+{1\over 2}\over x})^\rho
   \int_{l+{1\over 2}\over x}^\infty e^{-\pi u^2}u^{2-\rho}du\\
   &=-4\pi ({l+{1\over 2}\over x})^3\int_0^1
 u^{2-\rho}e^{-\pi({l+{1\over 2}\over x}u)^2}du.
      \end{align*}
      \end{thm}

   When $\Re\rho<1/2$, the following improvement of Theorem \ref{thm1.1}
   is also obtained in Section 3.

 \begin{thm} \label{thm1.3} If $\Re\rho<1/2$, then
  \[\|{1\over x}\sum_{m=1}^l\frak FZ^{-1}F_\rho({m\over x})-F_\rho(x)
 -\delta_{l,F_\rho}(x)\|_{L^2(0,\infty)}\ll {1\over 2l+1}\]
 for each positive integer $l$.
 \end{thm}

\section{Proof of Theorem \ref{thm1.2}}

 \begin{lem}\label{lem2.1} (\cite[Lemma 3.2]{li3})
  For every element $F\in\mathcal H_-$ we have
  \[\delta_{l,F}(x)=\int_0^{l+{1\over 2}\over x}\frak F Z^{-1}F(t)dt.\]
   \end{lem}

  \begin{lem}\label{lem2.2}  We can write
 \[\delta_{l,F_\rho}(x)
   ={\widehat\eta(1-\rho)(l+{1\over 2})^\rho
   \over \rho x^\rho}-\int_0^{l+{1\over 2}\over x}
  \left(\int_1^\infty \eta(tu)u^{-\rho}du\right)dt.\]
  \end{lem}

 \noindent\textit{Proof}. Since $\eta\in\mc H_\cap$, for any fixed $x>0$ and
 for all $t\geq 1$ we have $|\eta(tx)|\leq c_x (tx)^{-2}$
 for a constant $c_x$ depending on $x$. Thus,
 \[\int_1^\infty \left|\sum_{n\geq X+1}\eta(ntx)\right|t^{\Re\rho-1}dt
 \leq \frac{c_x}{x^2X(2-\Re\rho)}\to 0\]
as $X\to\infty$.  This implies that we can change the order of summation
and integration to write
 \[F_\rho(x)=Z\int_1^\infty \eta(tx)t^{\rho-1}dt.\]
   By Lemma \ref{lem2.1},
  \begin{align*}\delta_{l,F_\rho}(x)
  &=\int_0^\infty Z^{-1}F_\rho(t){\sin{l+{1\over 2}\over x}2\pi t \over\pi t}dt\\
  &=2\int_0^\infty dt\int_0^{l+{1\over 2}\over x}\cos 2\pi vtdv
     \int_1^\infty \eta(ut)u^{\rho-1}du\\
  &=2\int_0^{l+{1\over 2}\over x}dv\int_0^\infty \left(
     \int_1^\infty \eta(ut)u^{\rho-1}du\right)\cos 2\pi vtdt\\
  &=2\int_0^{l+{1\over 2}\over x}dv \int_1^\infty u^{\rho-1}du
  \int_0^\infty \eta(ut)\cos 2\pi vtdt\end{align*}
  where changes of order of integration are permissible because the
  triple integral is absolute integrable.
  Since $\frak F\eta=\eta$, we have
  \[2\int_0^\infty \eta(ut)\cos 2\pi vtdt=\frak F_t\{\eta(ut)\}(v)
  ={1\over u}\eta({v\over u})\]
  Thus
  \begin{align*}\delta_{l,F_\rho}(x) &=\int_0^{l+{1\over 2}\over x}dv
  \int_1^\infty \eta({v\over u})u^{\rho-2}du\\
  &=\int_0^{l+{1\over 2}\over x}v^{\rho-1}dv
  \int_0^\infty \eta({1\over u})u^{\rho-2}du-\int_0^{l+{1\over 2}\over x}dv
  \int_0^1 \eta({v\over u})u^{\rho-2}du\\
  &={\widehat\eta(1-\rho)(l+{1\over 2})^\rho\over\rho x^\rho}
  -\int_0^{l+{1\over 2}\over x}dt\int_1^\infty \eta(tu)u^{-\rho}du.
   \end{align*}

    This completes the proof of the lemma. $\hfill\Box$

  \noindent\textit{Proof of Theorem \ref{thm1.2}}.  Since
  \[{d\over dt}\left(-4\pi t^3e^{-\pi t^2}\right)=-4\pi(3t^2-2\pi t^4)
  e^{-\pi t^2}=8\pi t^2(\pi t^2-{3\over 2})e^{-\pi t^2}=\eta(t),\]
   we have
   \[\int_0^A\eta(t)dt=-4\pi{A^3\over e^{\pi A^2}}.\]
        Thus
  \begin{align*}&\int_0^{l+{1\over 2}\over x}
 \left(\int_1^\infty \eta(tu)u^{-\rho}du\right)dt
 =\int_1^\infty u^{-\rho}du\int_0^{l+{1\over 2}\over x}\eta(tu)dt\\
 &=\int_1^\infty u^{-\rho-1}du\int_0^{{l+{1\over 2}\over x}u}\eta(t)dt
  =-4\pi ({l+{1\over 2}\over x})^3\int_1^\infty
 u^{2-\rho}e^{-\pi({l+{1\over 2}\over x}u)^2}du\\
 &=-4\pi ({l+{1\over 2}\over x})^\rho\int_{l+{1\over 2}\over x}^\infty
   e^{-\pi u^2}u^{2-\rho}du\\
 &={\widehat\eta(1-\rho)(l+{1\over 2})^\rho\over\rho x^\rho}
 +4\pi ({l+{1\over 2}\over x})^\rho\int_0^{l+{1\over 2}\over x}
   e^{-\pi u^2}u^{2-\rho}du.
 \end{align*}
 By Lemma \ref{lem2.2},
 \[\delta_{l,F_\rho}(x)={\widehat\eta(1-\rho)(l+{1\over 2})^\rho
   \over \rho x^\rho}-\int_0^{l+{1\over 2}\over x}
  \left(\int_1^\infty \eta(tu)u^{-\rho}du\right)dt.\]
  It follows that
   \[\delta_{l,F_\rho}(x)= {\widehat\eta(1-\rho)(l+{1\over 2})^\rho
   \over \rho x^\rho}+4\pi ({l+{1\over 2}\over x})^\rho
   \int_{l+{1\over 2}\over x}^\infty e^{-\pi u^2}u^{2-\rho}du=
   -4\pi ({l+{1\over 2}\over x})^3\int_0^1
 u^{2-\rho}e^{-\pi({l+{1\over 2}\over x}u)^2}du.\]

      This completes the proof of Theorem \ref{thm1.2}.$\hfill\Box$

\section{Proof of Theorem \ref{thm1.3}}

 \begin{lem}\label{lem3.1} (\cite[Lemma 4.2]{li3})
  Let $f\in L^1(\mathbb R)$ be an even function.
 If $\frak Ff(0)=0$, then
 \[{1\over x}\sum_{m=1}^l\frak F  f({m\over x})
 =\sum_{m=1}^\infty {f(m x+0)+f(m x-0)\over 2}+\int_0^{1\over 2}
 f( x t){\sin(2l+1)\pi t\over\sin \pi t}dt+R(f, x)\]
   where
  \begin{align*}R(f,  x)&=\sum_{m=1}^\infty\int_m^{m+{1\over 2}}
 \left\{f(t x)-f(m x+0)\right\}{\sin(2l+1)\pi\over\sin \pi t}dt\\
 &+\sum_{m=1}^\infty\int_{m-{1\over 2}}^m
 \left\{f(t x)-f(m x-0)\right\}{\sin(2l+1)\pi t\over\sin \pi t}dt
 \end{align*}
 for every positive integer $l$.
 \end{lem}

  \begin{lem}\label{lem3.2} Let $f(t)=Z^{-1}F_\rho(t)$ with $\Re\rho<1/2$.
   Then
 \[\int_0^\infty |\sum_{m=1}^\infty
 \int_{-{1\over 2}}^{1\over 2}\left(f(t x+m x)-f(m x)\right)
 {\sin(2l+1)\pi t\over\sin \pi t}dt|^2dx\ll {1\over (2l+1)^4}\]
 \end{lem}

 \noindent\textit{Proof}.  By Mellin's inversion formula,
 \[ f(t)={1\over 2\pi i}\int_{c-i\infty}^{c+i\infty}
 {\widehat F_\rho(s)\over \zeta(s)}t^{-s}ds\]
 for $c\geq 1$.  We write
 \[(t x+m x)^{-s}-(m x)^{-s}
 =-s\int_{m x}^{t x+m x}{du\over u^{s+1}}
 ={-s\over (m x)^s}\int_1^{1+{t\over m}}{du\over u^{s+1}}\]
 Since
 \begin{align*}&|\left(f(t x+m x)-f(m x)\right)
 {\sin(2l+1)\pi t\over\sin \pi t}|\\
 &\leq |\left({1\over 2\pi i}\int_{c-i\infty}^{c+i\infty}
 {-s\widehat  F_\rho(s)\over \zeta(s)}
 [\int_{m x}^{t x+m x}u^{-s-1}du]ds\right)
 {\sin(2l+1)\pi t\over\sin \pi t}|\\
 &\leq{(2l+1)\pi t\over 2\pi\cdot 2t}{t\over m^{1+c} x^c}
 \int_{c-i\infty}^{c+i\infty}
 |{-s\widehat  F_\rho(s)\over \zeta(s)}||ds|\ll {1\over m^{1+c}}
 \end{align*}
 for $|t|\leq 1/2$ with the implied constant
 depending only on $ x$ and $l$ not on $t$, the series
  \[\sum_{m=1}^\infty\left(f(t x+m x)
 -f(m x)\right){\sin(2l+1)\pi t\over\sin \pi t}\]
  is uniformly convergent for $|t|\leq 1/2$.
  Thus
 \begin{align*}&\sum_{m=1}^\infty\int_{-{1\over 2}}^{1\over 2}\left(
 f(t x+m x)-f(m x)\right){\sin(2l+1)\pi t\over\sin \pi t}dt\\
 &={1\over 2\pi i}\int_{c-i\infty}^{c+i\infty}
 {-s\widehat  F_\rho(s)\over  x^s\zeta(s)}\{\sum_{m=1}^\infty{1\over m^s}
 \int_{-{1\over 2}}^{1\over 2}{\sin(2l+1)\pi t\over\sin \pi t}dt
 \int_1^{1+{t\over m}}{du\over u^{s+1}}\}ds\end{align*}
 where the change of order of integration is permissible as the triple
 integral is absolute integrable.

    By \eqref{eq1.1}, $\widehat F_\rho(s)=\xi(s)/(s-\rho)$
    for all complex $s$.  Hence,
  \[{\widehat F_\rho(s)\over\zeta(s)}={s(s-1)\pi^{-s/2}\Gamma(s/2)\over s-\rho}\]
  for all complex $s\neq\rho$.
     An extension of Stirling's
formula \cite[line 13, p. 151]{titch} is that
  for any fixed value of $\sigma$
 \[|\Gamma(\sigma+it)|\sim \sqrt{2\pi}e^{-{\pi|t|\over 2}}|t|^{\sigma-{1\over 2}}\]
  as $|t|\to\infty$.
  Thus,
  \[{\widehat F_\rho(s)\over\zeta(s)}\ll |s|^{-k}\]
   for any positive integer $k$ when $0\leq\Re s\leq 1$ and $|s|\to\infty$.
It follows that we can move the
 line of integration $\Re s=c\geq 1$ to $\Re s=\Re\rho+\epsilon<1/2$
 when $0< x<1$ and derive
 \begin{align*}&\sum_{m=1}^\infty\int_{-{1\over 2}}^{1\over 2}\left(
 f(t x+m x)-f(m x)\right){\sin(2l+1)\pi t\over\sin \pi t}dt\\
 &=\begin{cases} {1\over 2\pi i}\int_{\Re s=\Re\rho+\epsilon}
 {-s\widehat  F_\rho(s)\over  x^s\zeta(s)}\{\sum_{m=1}^\infty{1\over m^s}
 \int_{-{1\over 2}}^{1\over 2}{\sin(2l+1)\pi t\over\sin \pi t}dt
 \int_1^{1+{t\over m}}{du\over u^{s+1}}\}ds &\text{if $0< x<1$}\\
  {1\over 2\pi i}\int_{\Re s=1}
 {-s\widehat  F_\rho(s)\over  x^s\zeta(s)}\{\sum_{m=1}^\infty{1\over m^s}
 \int_{-{1\over 2}}^{1\over 2}{\sin(2l+1)\pi t\over\sin \pi t}dt
 \int_1^{1+{t\over m}}{du\over u^{s+1}}\}ds &\text{if $ x\geq 1$}.
 \end{cases} \end{align*}

   We write
 \begin{align*}&\int_{-{1\over 2}}^{1\over 2}{\sin(2l+1)\pi t\over\sin \pi t}dt
 \int_1^{1+{t\over m}}{du\over u^{s+1}}
 =\int_{-{1\over 2}}^{1\over 2}{\cos(2l+1)\pi t\over s(2l+1)\pi}
 {d\over dt}\{{{1-(1+{t\over m})^{-s}\over t}\over {\sin\pi t\over t}}\}dt\\
 &=\int_{-{1\over 2}}^{1\over 2}{\cos(2l+1)\pi t\over(2l+1)\pi}
 \{{{t\over m}(1+{t\over m})^{-s-1}-{1-(1+{t\over m})^{-s}\over s}
 \over t\sin\pi t}+{d\over dt}({t\over\sin\pi t})
 {1\over t}\int_1^{1+{t\over m}}{du\over u^{s+1}}\}dt\\
 &=\int_{-{1\over 2}}^{1\over 2}{\cos(2l+1)\pi t\over(2l+1)\pi}
 \{{{(1+{t\over m})^{-s-1}\over t}\int_1^{1+{t\over m}}
 {1-({1+{t\over m}\over u})^{s+1}\over t}du
 \over {\sin\pi t\over t}}+{d\over dt}({t\over\sin\pi t})
 {1\over t}\int_1^{1+{t\over m}}{du\over u^{s+1}}\}dt\\
 &=\int_{-{1\over 2}}^{1\over 2}{\cos(2l+1)\pi t\over(2l+1)\pi}
 \{{{s+1\over t^2(1+{t\over m})^{s+1}}\int_1^{1+{t\over m}}du
 \int_1^{1+{t\over m}\over u}{dv\over v^{s+2}}\over {\sin\pi t\over t}}
 +{d\over dt}({t\over\sin\pi t})
 {1\over t}\int_1^{1+{t\over m}}{du\over u^{s+1}}\}dt\\
 &=\int_{-{1\over 2}}^{1\over 2}{\cos(2l+1)\pi t\over(2l+1)\pi}
 \{{{s+1\over m^2(1+{t\over m})^{s+1}}\int_0^1du
 \int_0^{1-u\over 1+{t\over m}u}{dv\over (1+{t\over m}v)^{s+2}}\over {\sin\pi t\over t}}
 +{d\over dt}({t\over\sin\pi t})
 {1\over m}\int_0^1{du\over (1+{t\over m}u)^{s+1}}\}dt
  \end{align*}
  Since $\sin t\geq 2t/\pi$ for $t\in [0,\pi/])$, for $|t|\leq 1/2$
   we have $|\sin\pi t/t|\geq 2$.  Applying partial integration to the last term
   of the above identity we obtain that
  \[\int_{-{1\over 2}}^{1\over 2}{\sin(2l+1)\pi t\over\sin \pi t}dt
 \int_1^{1+{t\over m}}{du\over u^{s+1}}\ll {(|s|+1)^2\over m(2l+1)^2}.\]
  It follows that
  \[|\sum_{m=1}^\infty\int_{-{1\over 2}}^{1\over 2}\left(
 f(t x+m x)-f(m x)\right){\sin(2l+1)\pi t\over\sin \pi t}dt|
 \ll \begin{cases} { x^{-\Re\rho-\epsilon}\over (2l+1)^2} &\text{if $0< x<1$}\\
  { x^{-1}\over (2l+1)^2} &\text{if $ x\geq 1$}.\end{cases}\]
   Therefore
 \[\int_0^\infty |\sum_{m=1}^\infty
 \int_{-{1\over 2}}^{1\over 2}\left(f(t x+m x)-f(m x)\right)
 {\sin(2l+1)\pi t\over\sin \pi t}dt|^2dx\ll {1\over (2l+1)^4}.\]

  This completes the proof of the lemma.
 $\hfill\Box$

  \begin{lem}\label{lem3.3}   Let $f(t)=Z^{-1}F_\rho(t)$.  Then
  \[\int_0^\infty |\int_0^{1\over 2} f(xt)({1\over\pi t}-
 {1\over\sin\pi t})\sin(2l+1)\pi tdt|^2dx\ll{1\over (2l+1)^2}.\]
 \end{lem}

 \noindent\textit{Proof}.  Note that
 \[f(t)={1\over 2\pi i}\int_{\Re s=c}{\widehat F_\rho(s)\over \zeta(s)}t^{-s}ds\]
 for $c\geq 1$.  Because the following double integral is absolute integrable,
 we can change order of integration to derive
 \begin{align*}&\int_0^{1\over 2} f(xt)({1\over\pi t}-
 {1\over\sin\pi t})\sin(2l+1)\pi tdt\\
 &={1\over 2\pi i}\int_{\Re s=c}{\widehat F_\rho(s)\over \zeta(s)}\{
 \int_0^{1\over 2} t^{-s}({1\over\pi t}-{1\over\sin\pi t})
 \sin(2l+1)\pi tdt\}{ds\over x^s}.
 \end{align*}
 Similarly as in proof of Lemma \ref{lem3.2} we write
 \begin{align*}&\int_0^{1\over 2} f(xt)({1\over\pi t}-
 {1\over\sin\pi t})\sin(2l+1)\pi tdt\\
 &=\begin{cases} {1\over 2\pi i}\int_{\Re s=\Re\rho+\epsilon}
 {\widehat F_\rho(s)\over \zeta(s)}\{
 \int_0^{1\over 2} t^{-s}({1\over\pi t}-{1\over\sin\pi t})
 \sin(2l+1)\pi tdt\}{ds\over x^s} &\text{if $0<x<1$}\\
 {1\over 2\pi i}\int_{\Re s={1\over 2}+\epsilon}
 {\widehat F_\rho(s)\over \zeta(s)}\{
 \int_0^{1\over 2} t^{-s}({1\over\pi t}-{1\over\sin\pi t})
 \sin(2l+1)\pi tdt\}{ds\over x^s} &\text{if $x\geq 1$}\end{cases}
 \end{align*}
 for a sufficiently small positive $\epsilon$.

 By \cite[(13), p. 190]{ahl},
 \[{1\over t}-{\pi\over\sin(\pi t)}=2t\sum_{n=1}^\infty {(-1)^n\over n^2-t^2}
 =2t\sum_{n=1}^\infty{(-1)^n\over n^2}\sum_{k=0}^\infty ({t\over n})^{2k}\]
 for $0<t<1/2$.   Note that
 \[2\sum_{n=1}^\infty {(-1)^n\over n^2}=-{\pi^2\over 6}\]
 Applying partial integration twice we derive that
   \begin{align*}&\int_0^{1\over 2} t^{-s}\{2t\sum_{n=1}^\infty{(-1)^n\over n^2}
   \sum_{k=1}^\infty ({t\over n})^{2k}\}\sin(2l+1)\pi tdt\\
   &=-\int_0^{1\over 2} t^{3-s}\{2\sum_{n=1}^\infty{(-1)^n\over n^2}
   \sum_{k=1}^\infty ({t\over n})^{2k-2}\}d{\cos(2l+1)\pi t\over (2l+1)\pi}
   \ll {(|s|+3)^2\over (2l+1)^2}
   \end{align*}
 for $0<\Re s<1$.

  Also, by partial integration we have
 \begin{align*}
 &\int_0^{1\over 2} t^{-s}\{2t\sum_{n=1}^\infty {(-1)^n\over n^2}\}\sin(2l+1)\pi t\,dt
 =-{\pi^2\over 6}\int_0^{1\over 2} t^{1-s}\sin(2l+1)\pi tdt\\
 &={\pi(s-1)\over 6(2l+1)}\int_0^{1\over 2} t^{-s}\cos(2l+1)\pi tdt
 ={\pi^s(s-1)\over 6(2l+1)^{2-s}}\int_0^{l+{1\over 2}} t^{-s}\cos t\,dt.
 \end{align*}
 for $0<\Re s<1$.   By \cite[Example 10, p. 162]{titch},
 \[\int_0^\infty t^{-s}\cos t\,dt=\Gamma(1-s)\sin{\pi s\over 2}\]
 for $0<\Re s<1$.  As
 \[\int_{l+{1\over 2}}^\infty t^{-s}\cos t\,dt=-{\sin(l+{1\over 2})\over (l+{1\over 2})^s}
 +s\int_{l+{1\over 2}}^\infty t^{-s-1}\sin t\,dt\ll {|s|+1\over (l+{1\over 2})^{\Re s}}\]
 we have
 \[\int_0^{1\over 2} t^{-s}\{2t\sum_{n=1}^\infty {(-1)^n\over n^2}\}\sin(2l+1)\pi t\,dt
 =-{\pi^s\Gamma(2-s)\sin{\pi s\over 2}\over 6(2l+1)^{2-s}}
 +O\left({(|s|+1)^2\over (2l+1)^2}\right)\]
 Thus, we have obtained that
 \[\int_0^{1\over 2} t^{-s}({1\over\pi t}-{1\over\sin\pi t})
 \sin(2l+1)\pi tdt =-{\pi^s\Gamma(2-s)\sin{\pi s\over 2}\over 6(2l+1)^{2-s}}
 +O\left({(|s|+3)^2\over (2l+1)^2}\right)\]
 Hence
 \begin{align*}&\int_0^{1\over 2} f(xt)({1\over\pi t}-
 {1\over\sin\pi t})\sin(2l+1)\pi tdt\\
 &={1\over 2\pi i}\int_{\Re s=c}{\widehat F_\rho(s)\over \zeta(s)}\{
 -{\pi^s\Gamma(2-s)\sin{\pi s\over 2}\over 6(2l+1)^{2-s}}
 +O\left({(|s|+3)^2\over (2l+1)^2}\right)\}{ds\over x^s}.
 \end{align*}

 By the proof of Lemma \ref{lem3.2},
  \[{\widehat F_\rho(s)\over\zeta(s)}\ll |s|^{-k}\]
  for any positive integer $k$ when $0\leq\Re s\leq 1$ and $|s|\to\infty$.
 Thus, by estimations
 \[{1\over 2\pi i}\int_{\Re s=c}{\widehat F_\rho(s)\over \zeta(s)}\{
 -{\pi^s\Gamma(2-s)\sin{\pi s\over 2}\over 6(2l+1)^{2-s}}
 +O\left({(|s|+3)^2\over (2l+1)^2}\right)\}{ds\over x^s}
 \ll \begin{cases} {x^{-\Re\rho-\epsilon}\over (2l+1)^{2-\Re\rho-\epsilon}}
 &\text{if $0<x<1$}\\
  {x^{-1}\over 2l+1} &\text{if $x\geq 1$}.\end{cases}\]
 Therefore,
  \[\int_0^{1\over 2} f(xt)
  ({1\over\pi t}-{1\over\sin\pi t})\sin(2l+1)\pi tdt
   =O\left(\begin{cases} {x^{-\Re\rho-\epsilon}\over (2l+1)^{2-\Re\rho-\epsilon}}
   &\text{if $0<x<1$}\\
  {x^{-1}\over 2l+1} &\text{if $x\geq 1$}\end{cases}\right).\]
 It follows that
  \[\int_0^\infty |\int_0^{1\over 2} f(xt)({1\over\pi t}-
 {1\over\sin\pi t})\sin(2l+1)\pi tdt|^2dx\ll{1\over (2l+1)^2}.\]

 This completes the proof of the lemma.
 $\hfill\Box$

   \begin{lem}\label{lem3.4} Let $f(t)=Z^{-1}F_\rho(t)$ with
   $\Re\rho<1/2$.  Then
  \[\int_0^\infty\left|\int_{1\over 2}^\infty
 f(xt){\sin(2l+1)\pi t\over\pi t}dt\right|^2dx\ll {1\over (2l+1)^4}\]
 \end{lem}

 \noindent\textit{Proof}.  As
 \[f(t)={1\over 2\pi i}\int_{\Re s=c}{\widehat F_\rho(s)\over \zeta(s)}t^{-s}ds\]
 for $c\geq 1$, we can write
 \[\int_{1\over 2}^\infty f(xt){\sin(2l+1)\pi t\over\pi t}dt
 ={1\over 2\pi i}\int_{\Re s=c}{\widehat F_\rho(s)\over \zeta(s)}
 \{\int_{1\over 2}^\infty{\sin(2l+1)\pi t\over\pi t^{1+s}}dt\}{ds\over x^s}\]
 where change of order of integration is permissible as the double integral is
 absolute integrable.

 By partial integration,
 \[\int_{1\over 2}^\infty f(xt){\sin(2l+1)\pi t\over\pi t}dt
 =\begin{cases} {-1\over 2\pi i}\int_{\Re s=1}{(1+s)\widehat F_\rho(s)\over \zeta(s)}
 \{\int_{1\over 2}^\infty{\cos(2l+1)\pi t\over(2l+1)\pi^2 t^{2+s}}dt\}{ds\over x^s}
 &\text{if $x\geq 1$}\\
  {-1\over 2\pi i}\int_{\Re s=\Re\rho+\epsilon}{(1+s)\widehat F_\rho(s)\over \zeta(s)}
 \{\int_{1\over 2}^\infty{\cos(2l+1)\pi t\over(2l+1)\pi^2 t^{2+s}}dt\}{ds\over x^s}
 &\text{if $x<1$}\end{cases}\]
 where $\epsilon$ is a sufficiently small positive satisfying
 $\Re\rho+\epsilon<{1\over 2}$.   Also,
 \[\int_{1\over 2}^\infty{\cos(2l+1)\pi t\over(2l+1)t^{2+s}}dt
 =-{2^{2+s}\sin(l+{1\over 2})\pi\over (2l+1)^2\pi}+{2+s\over \pi(2l+1)^2}
 \int_{1\over 2}^\infty {\sin(2l+1)\pi t\over t^{3+s}}dt\ll {|s|+2\over (2l+1)^2}\]
 for $0<\Re s<1$.  It follows that
 \[|\int_{1\over 2}^\infty f(xt){\sin(2l+1)\pi t\over\pi t}dt|
 \ll \begin{cases}{x^{-\Re\rho-\epsilon}\over (2l+1)^2}&\text{if $0<x\leq 1$}\\
 {x^{-1}\over (2l+1)^2}&\text{if $x>1$}.\end{cases}\]
 This inequality implies that
 \[\int_0^\infty\left|\int_{1\over 2}^\infty
 f(xt){\sin(2l+1)\pi t\over\pi t}dt\right|^2dx\ll {1\over (2l+1)^4}.\]

  This completes the proof of the lemma.
 $\hfill\Box$

 \noindent\textit{Proof of Theorem \ref{thm1.3}}. For $f(t)=Z^{-1}F_\rho(t)$,
 we have
 $f(0)=0$ and $\frak F  f(0)=0$.
 Also, $f$ is an even function on $\mathbb R$. Thus, by Lemma \ref{lem3.1}
 \begin{equation}\label{eq3.4}{1\over x}\sum_{m=1}^l\frak F  f({m\over x})
 =\sum_{m=1}^\infty {f(mx+0)+f(mx-0)\over 2}
 +\int_0^{1\over 2}f(x t){\sin(2l+1)\pi t\over\sin \pi t}dt+R(f,x)
 \end{equation}
  where
  \[R(f, x)=\sum_{m=1}^\infty\int_{-{1\over 2}}^{1\over 2}
 \left(f(tx+mx)-f(mx)\right){\sin(2l+1)\pi t\over\sin \pi t}dt.\]
 We rewrite \eqref{eq3.4} as
  \begin{align}\label{eq3.5}
   &{1\over x}\sum_{m=1}^l\frak FZ^{-1}F_\rho({m\over x})-\delta_{l,F_\rho}(x)-F_\rho(x)\notag\\
   &=\int_0^{1\over 2} f(xt)({1\over\sin\pi t}-{1\over\pi t})\sin(2l+1)\pi tdt
   -\int_{1\over 2}^\infty f(xt){\sin(2l+1)\pi t\over\pi t}dt\\
 &+\sum_{m=1}^\infty\int_{-{1\over 2}}^{1\over 2}\left(f(t x+m x)-f(m x)\right)
 {\sin(2l+1)\pi t\over\sin \pi t}dt.\notag
  \end{align}

   By \eqref{eq3.5}, the Minkowski's inequality \cite[Theorem 3.5, p. 63]{rudin}, and
    Lemmas \ref{lem3.2}, \ref{lem3.3} and \ref{lem3.4} we obtain that
   \[\|({1\over x}\sum_{m=1}^l\frak FZ^{-1}F_\rho({m\over x})
 -\delta_{l,F_\rho}(x))-F_\rho(x)\|_{L^2(0,\infty)}\ll {1\over 2l+1}.\]

    This completes the proof of Theorem \ref{thm1.3}. $\hfill\Box$

\end{document}